\documentclass[11pt]{article}
\usepackage[utf8]{inputenc}
\usepackage{fullpage}
\usepackage{amsmath,amssymb,amsthm,geometry}
\usepackage{xcolor}
\usepackage{tikz,tikz-3dplot}
\usepackage[colorlinks=true,linkcolor=blue]{hyperref}
\hypersetup{
   citecolor=blue,
   linkcolor=blue,
   }
\usepackage{hyperref}
\usepackage[psamsfonts]{eucal}
\usepackage{microtype}
\usepackage{color} 
\usepackage{mathtools}

\newcommand{\R}{\mathbb{R}}

 \usepackage{url}
 \usepackage{framed}
\usepackage{latexsym}
 \usepackage{graphicx, epstopdf, verbatim}
\usepackage{psfrag}

\newtheorem{theorem}{Theorem}[section]

\newtheorem{lemma}[theorem]{Lemma}

\theoremstyle{definition}\newtheorem{assertion}{Assertion}

\usepackage{mathtools}
\usepackage{setspace}

\usepackage{breqn}
\usepackage{wasysym}

\usepackage{tikz}

\usepackage{enumerate} 

\usepackage{empheq}

\title{Invariant Curves for Degenerate Hyperbolic Maps of the Plane}
\usepackage{times, verbatim}
\author{Charles Fefferman\footnote{Princeton University, Mathematics Department, Fine Hall, Washington Road, Princeton NJ, 08544-1000, USA. Partially supported by NSF grant DMS-1700180}} 
\date{September 22, 2021}
\begin{document}
\maketitle
\begin{abstract}We prove existence and uniqueness of an unstable manifold  
for a degenerate hyperbolic map of the plane arising in statistics.
\end{abstract}
\begin{center}
\textbf{AMS Classification 37E30}
\end{center}

\textbf{Note:} After posting an earlier version of this paper I learned that  
the results presented here are special cases of theorems due to I.  
Baldoma, E. Fontich, R. de la Llave, and M. Pau. An appendix to this  
paper provides details. I am grateful to Baldoma, Fontich and Pau for  
supplying the appendix, and to D. Cordoba and R. de la Llave for  
making me aware of the prior results.

\section{Introduction}\label{sec:intro}
The standard stable manifold theorem \cite{Hirsch1977} applies in particular to maps of the plane having the form 
\[ \Phi:(x,y) \mapsto (X,Y) \]
with 

\begin{align}
\label{eqn:m1}
X &= \lambda_1 x + O(|(x,y)|^2) \\ 
Y &= \lambda_2 y + O(|(x,y)|^2) \nonumber
\\ \nonumber
&|\lambda_1| > 1 > |\lambda_2|. 
\end{align}

\noindent
Such a map has an invariant unstable manifold tangent to the x-axis.
Here we study a degenerate case in which $\lambda_1 = + 1$ and $\lambda_2 = -1$. More precisely, we study smooth maps $\Phi:(x,y) \mapsto (X,Y)$, with 
\begin{align}
\label{eqn:m2} 
X &= x + x^2+\mu xy + O(|(x,y)|^3)\\ \nonumber
Y &= - y + \lambda xy + O(|(x,y)|^3)\\ \nonumber
&\mu \in \mathbb{R}, \lambda > 0. 
\end{align}

\noindent
Like \eqref{eqn:m1}, a map \eqref{eqn:m2} expands $x$ and contracts $y$ in the region of interest (where, in particular, $x>0$ and $|y|<<x$), but the stretching and shrinking arise from second order terms in the Taylor expansion of $\Phi$. We will prove existence and uniqueness of a smooth invariant curve, tangent to the positive x-axis, for maps of the form (2). Note that $\Phi^{-1}$ is defined in a neighborhood of the origin, since $\Phi'(0,0)$ is invertible. 


%
%
%
%

Our interest in maps \eqref{eqn:m2} arises from Lee et al. \cite{lee2020valid}, which proposes critical values for the $t$-ratio associated with the method of instrumental variables regression, which has received a great deal of attention in economics research, and has been widely employed across many empirical disciplines. In the instrumental variable model, there is an outcome of interest $Y$ (e.g., rates of illness), a causal factor of interest $X$ (e.g., receipt of a vaccine) and an “instrument” $Z$ (e.g., random assignment to either receiving or not receiving encouragement to take the vaccine). Under certain conditions, the causal effect of interest (e.g., the impact of vaccine receipt on rates of illness) has been shown to be equal to the ratio of two regression coefficients: the coefficient in a regression of $Y$ on $Z$ divided by that from the regression of $X$ on $Z$. Empirical researchers have typically used a $t$-ratio to test statistical significance using a constant critical value threshold (e.g. $\pm1.96$ for a 5\% test), presuming that the $t$-ratio is approximately standard normal. But as discussed in \cite{lee2020valid}, that approximation has been shown in the economics literature to be quite poor in many empirically-relevant cases. \cite{lee2020valid} corrects for this by providing critical values that depend on the $F$-statistic associated with the regression of $X$ on $Z$. A piece of this critical value function is the solution to a functional equation, which takes the form \eqref{eqn:m2}, after a change of variables described in \cite{lee2020valid}. Our main result on maps \eqref{eqn:m2} provides a rigorous proof of the existence and uniqueness of such a solution to the functional equation, the numerical solution of which was previously reported in an earlier version of \cite{lee2020valid}.

\bigskip
\noindent
The precise statement of our result is as follows. Here and below, "smooth" means $C^{\infty}$.

\begin{theorem}\label{thm1.1}
Let $\Phi: U \mapsto \mathbb{R}^2$ be a smooth map defined in a neighborhood U of the origin in $\mathbb{R}^2$. Suppose $\Phi$ has the form (2). Then there exists a smooth curve 
\[\Gamma = \{(x,F(x)): x \in [0,\delta] \} \subset U \] with the following properties. 
\begin{enumerate}[(I)]
    \item Invariance: $\Phi^{-1}(\Gamma) \subset \Gamma$.
    \item Tangency: $F(x)=O(x^3)$ as $x \rightarrow 0$. 
    \item Uniqueness: Let $\tilde{\Gamma} = \{ (x,\tilde{F}(x)): x \in [0,\tilde{\delta}] \}$, where $\Phi^{-1}(\tilde{\Gamma}) \subset \tilde{\Gamma}$ and $x^{-2/3}\tilde{F}(x) \rightarrow 0$ as $x \rightarrow 0^{+}$. Then $\tilde{F}=F$ on $[0,\delta^{\#}]$ for some $\delta^{\#}>0$.
\end{enumerate}

\end{theorem}

\noindent
We now sketch the proof of the above Theorem. We first make a change of variable to bring $\Phi$ to the form $\Phi:(x,y) \mapsto (X,Y)$ with 

\begin{align}
X &= x + x^2 + x^3 \theta_A(x,y) + y\theta_B(x,y) , 
\label{eqn:m3} 
\\
\label{eqn:m4}
Y &= -y(1 - \lambda x+ x^2 \theta_C(x,y) + y \theta_D(x,y)) + x^{N+100}\theta_E(x,y) .
\end{align}

\noindent
Here the $\theta$'s are smooth functions, and N is as large as we please. But for the term $x^{N+100}\theta_E(x,y)$ in \eqref{eqn:m4}, the x-axis would be an invariant curve for $\Phi$. We look for an invariant curve of the form 
\[\Gamma = \{ (x,F(x)): x \in [0,\delta])\}\] with 
\begin{align} \label{eqn:m5}
\left|\left(\frac{d}{dx}\right)^m F(x)\right| \leq K_m x^{N-m} 
\end{align}
on $[0,\delta]$, $m=0,1,\ldots, N-10$,  for carefully selected constants 
 $K_0, K_1, ..., K_{N-10}$.

 To produce such an F, we proceed as follows. Fix a small number $\rho$, $0 < \rho << \delta$. Later, we will let $\rho$ tend to zero. We start with a horizontal line segment 
 
 \[\Gamma_{\rho}^{0} = \{(x,0): x \in [0,\rho]\}\]

\noindent
and then apply an iterate of $\Phi$ to produce the image 
\[\Gamma_{\rho} = \Phi^{\bar{\nu}}(\Gamma_{\rho}^{0})\]
where we take $\bar{\nu}$ to be the least integer for which $\Phi^{\bar{\nu}}(\Gamma_{\rho}^{0})$ contains points (x,y) with $x > \delta$. Note that $\bar{\nu} \rightarrow \infty$ as $ \rho \rightarrow 0$. We will show that $\Gamma_{\rho}$ has the form:
\begin{align} 
 \Gamma_{\rho} = \{ (x, F_{\rho}(x)): x \in \left[0, x_{\rho}^{MAX}\right] \}  
\end{align} 
 with $x_{\rho}^{MAX} \geq \delta$ and 
 \begin{align} 
 \label{eqn:m7}
\left|\left(\frac{d}{dx}\right)^m F_{\rho}(x)\right| \leq K_m x^{N-m} \   \text{ on } [0,x_{\rho}^{MAX}], \ m=0,1,...,N-10. 
\end{align}

\noindent
Here, the $K_m$ are as in \eqref{eqn:m5}; they are independent of $\rho$. Moreover, we will show that $\Gamma_{\rho}$ is approximately invariant, in the sense that every point of $\Phi^{-1}(\Gamma_{\rho})$ lies within a distance $\rho$ of a point of $\Gamma_{\rho}$. By Ascoli's theorem, we can find a sequence $\rho_1, \rho_2, ...$ tending to zero, such that the curves $\Gamma_{\rho_{i}}$ tend to a limiting curve $\Gamma$ in the $C^{N-11}$ topology. That curve is invariant, highly tangent to the x-axis and $C^{N-11}$ for N as large as we please. 

The uniqueness and $C^{\infty}$ smoothness assertions of our theorem then follow easily. 

We now delve slightly deeper by providing a few words about the proof of \eqref{eqn:m7}, and that of the approximate invariance of the curves $\Gamma_{\rho}$. We will successively pick constants 
\[1 = K_0 << K_1 << ... << K_{N-10}\] and then pick a small enough $\delta$ depending on the $K_m$. For these constants, we study smooth curves of the form 
 \begin{align}\label{equ:8}
 \Gamma &= \{ (x,f(x)): x \in [0,x_{MAX}]\} \\
  \text{such that }
  \left|\left(\frac{d}{dx}\right)^m f(x) \right| &\leq K_m x^{N-m} \ \text{ on } [0, x_{MAX}],\  m = 0, 1, ..., N-10 .
  \nonumber 
 \end{align}

\noindent
We show that if $\Gamma $ is of the form \eqref{equ:8} and $x_{MAX} \leq \delta$, then $\Phi(\Gamma)$, the image of $\Gamma$ under $\Phi$, is again of the form \eqref{equ:8}, with different $x_{MAX}$ and f, but with the same $K_0,...,K_{N-10}$. Starting from the horizontal line segment $\Gamma_{\rho}^0$, which clearly has the form \eqref{equ:8}, we can therefore repeatedly take the image under $\Phi$, always preserving \eqref{equ:8}, until at last $x_{MAX}$ in \eqref{equ:8} exceeds $\delta$. Thus, we conclude that $\Gamma_{\rho} = \Phi^{\bar{\nu}}(\Gamma_{\rho}^0)$ is of the form \eqref{equ:8}. That's the plan of our proof of \eqref{eqn:m7}. 

To establish the approximate $\Phi$-invariance of our curve $\Gamma_{\rho}$, our main tool is the following Shadowing Lemma.  

\begin{lemma}  ~\\ 
Let $(x,y)$, $(\hat{x}, \hat{y}) \in \mathbb{R}^2$, 
with $0<x< \delta$, $|y| \leq K_0 x^N$, and 
$\frac{|x-\hat{x}|+ x^{-3}|y-\hat{y}|}{x^8} \leq 1$.  
\\ 
Then the points $(X,Y) = \Phi(x,y)$ and $(\hat{X}, \hat{Y}) = \Phi(\hat{x}, \hat{y})$ satisfy 
\[ 
\frac{|X-\hat{X}|+ X^{-3}|Y-\hat{Y}|}{X^8} \leq \frac{|x-\hat{x}| +x^{-3}|y-\hat{y}|}{x^8} \leq 1.
\] 
\label{lemma:1.2}
\end{lemma}

\noindent
The shadowing lemma allows us  to prove the approximate $\Phi$-invariance of $\Gamma_{\rho} $ by the following argument.
Let (x,y) $\in \Gamma_{\rho}$. By definition, $(x,y) \in \Phi^{\bar{\nu}}(\Gamma_{\rho}^0)$, i.e, (x,y) = $\Phi^{\bar{\nu}}(x_0,0) $ for some $x_0 \in [0, \rho]$. Let $(x_{\nu}, y_{\nu}) = \Phi^{\nu}(x_0, 0)$ for $\nu = 0, 1, ..., \bar{\nu}$. Note that all the $y_{\nu}$ satisfy $|y_{\nu}| \leq K_0 x_{\nu}^N$, thanks to the invariance of (8) under the map $\Phi$. We can easily find a point $(\tilde{x}_0, 0) \in \Gamma_{\rho}^0$ such that 
$(\hat{x}_0, \hat{y}_0) \equiv \Phi(\tilde{x}_0, 0)$ satisfies $\hat{x}_0 = x_0$, $|\hat{y}_0| \leq K_0 \hat{x}_0^N$, hence, 

\begin{align}\label{equ:9}
 \frac{|x_0-\hat{x}_0|+ x_0^{-3}|y_0-\hat{y}_0|}{x_0^8}  \leq K_0 x_0^{N-11} \leq K_0 \rho^{N-11} .
\end{align}

\noindent
Let $(\hat{x}_{\nu},\hat{y}_{\nu}) = \Phi^{\nu}(\hat{x}_{0}, \hat{y}_{0})$ for $\nu = 0, 1, ..., \bar{\nu}$. Starting from \eqref{equ:9}, and repeatedly applying the shadowing lemma, we learn that 

\[\frac{|x_{\bar{\nu}}-\hat{x}_{\bar{\nu}}|+ x_{\bar{\nu}}^{-3}|y_{\bar{\nu}}-\hat{y}_{\bar{\nu}}|}{x_{\bar{\nu}}^8}  \leq K_0 \rho^{N-11}. \]
In particular,
\begin{align}
\label{eqn:m10} 
|(x_{\bar{\nu}},y_{\bar{\nu}}) - (\hat{x}_{\bar{\nu}} , \hat{y}_{\bar{\nu}})| \leq \rho^{N-11}. 
\end{align}

\noindent
We now recall that 
\[(x_{\bar{\nu}},y_{\bar{\nu}}) = \Phi^{\bar{\nu}}(x_0, 0) = (x,y)\]
and that 
\[(\hat{x}_{\bar{\nu}},\hat{y}_{\bar{\nu}}) = \Phi^{\bar{\nu}}(\hat{x}_0, \hat{y}_0) = \Phi^{\bar{\nu}}(\Phi(\tilde{x}_0,0)) = \Phi(\Phi^{\bar{\nu}}(\tilde{x}_0,0)). \]

\noindent
Letting $(x^{\#},y^{\#}) = \Phi^{\bar{\nu}}(\tilde{x}_0,0) \in \Phi^{\bar{\nu}}(\Gamma_{\rho}^{0})=\Gamma_{\rho}$, we see that $(\hat{x}_{\bar{\nu}},\hat{y}_{\bar{\nu}})$ = $\Phi(x^{\#},y^{\#})$. 
Thus, \eqref{eqn:m10} shows that (x,y) lies within distance $\rho^{N-11}$ of a point $\Phi(x^{\#}, y^{\#})$, with $(x^{\#}, y^{\#}) \in \Gamma_{\rho}$. Since (x,y) here is an arbitrary point of $\Gamma_{\rho}$, this concludes the proof of approximate $\Phi$-invariance of $\Gamma_{\rho}$. 
This also concludes our summary of the proof of our theorem. 

In the sections below, we provide full details of the proof. We warn the reader that our notation in this introduction is not entirely consistent with the notation in subsequent sections. However, we have accurately summarized the main ideas. 

Dynamical systems researchers are aware that good things arising from the linearization of a map may also arise from higher terms in its Taylor expansion; Theorem \ref{thm1.1} is a case in point.

I thank Lai-Sang Young and Rafael de la Llave for useful comments on invariant manifolds. I'm grateful to the authors of Lee et al. \cite{lee2020valid} for posing an intriguing problem with a practical application, and to Peter Ozsvath for putting me in touch with them.

\section{Proof}
\subsection{Change of Coordinates}

\medskip

We look at maps

\begin{equation}
\label{eqn:I11} 
\Phi:(x,y) \longmapsto  (X,Y) 
\end{equation}
where
\begin{align}
\label{eqn:I12}
    Y &= -y(1-\lambda x + y \theta_{1} + x^{2}\theta_{2}) + x^{N}\theta_{3} ,
\\ 
\label{eqn:I13}
    X &= x + x^{2} + x^{3}\theta_{4} + y\theta_{5}.
\end{align}

\noindent
Here and below, $\theta$'s denote smooth functions of $(x, y), \lambda >0,$ and $N \geq 3$.
\\[1ex] 
Note that a map of the form \eqref{eqn:m2} has the form \eqref{eqn:I11}, \eqref{eqn:I12}, \eqref{eqn:I13} with $N=3$, since the $O(|(x,y)|^3)$ terms in \eqref{eqn:m2} may be expressed as $\theta_A(x,y)x^3+\theta_B(x,y)x^2y+\theta_C(x,y)xy^2+\theta_D(x,y)y^3$.
\\[1ex]  

\noindent
We make a change of variables: 
\begin{align*}
    \tilde{Y} = Y + \gamma X^{N}, \qquad \tilde{y} = y + \gamma x^{N}. \qquad (\gamma	\in \mathbb{R})
\end{align*}

\noindent
This changes $\Phi$  to a map
\begin{align*}
    \tilde{\Phi} : (x, \tilde{y}) \longmapsto (X, \tilde{Y} ).
\end{align*}

\noindent
We will show that, by picking the correct $\gamma$, we can arrange that $\tilde{\Phi}$ has the same form as $\Phi$ but with $N+1$ in place of $N$. We write $\theta_{i}$ for integers $i$ to denote smooth functions of $(x,y)$, or equivalently, smooth functions of $(x, \tilde{y})$. Our $\theta_{i}$ will be independent of $\gamma$.

\noindent
Note, $y = \tilde{y} - \gamma x^{N}$, so 
\begin{align*}
    Y = - [\tilde{y} - \gamma x^{N}] (1-\lambda x + [\tilde{y} - \gamma x^{N}]\theta_{1}+ x^{2}\theta_{2}) + x^{N}\theta_{3},
\end{align*}
where now the $\theta$'s are regarded as smooth functions of $(x, \tilde{y})$. So,
\begin{align*}
    Y &= -\tilde{y}(1-\lambda x + \tilde{y}\theta_{1} + x^{2}(\theta_{2}-\gamma x^{N-2}\theta_{1})) + \gamma x^{N}  \\ & \qquad + \gamma x^{N} (-\lambda x + [\tilde{y} - \gamma x^{N}]\theta_{1} + x^{2}\theta_{2}) + x^{N}\theta_{3}
    \\
     &= -\tilde{y}(1 - \lambda x + \tilde{y}\theta_{1} + x^{2}(\theta_{2}-\gamma x^{N-2}\theta_{1} - \gamma x^{N-2}\theta_{1})) + \gamma x^{N} + \theta_{3}x^{N}  \\ & \qquad + x^{N+1}(-\lambda \gamma -\gamma^{2}x^{N-1}\theta_{1} + x\gamma \theta_{2})
     \\
     &= -\tilde{y}(1-\lambda x + \tilde{y}\theta_{1} + x^{2}\theta_{6}) + (\gamma + \theta_{3})x^{N} + \theta_{7}x^{N+1} \text{ for smooth functions $\theta_{6}, \theta_{7}$.}
\\[1ex] 
\text{Similarly,}
\\[1ex] 
    \gamma X^{N} &= \gamma(x + x^2 + x^{3}\theta_{4} + y\theta_{5})^{N}
    \\
    &= \gamma(x + x^{2} + x^{3}\theta_{4} + [\tilde{y} - \gamma x^{N}]\theta_{5})^{N}
    \\
    &= \gamma(x + x^{2} + x^{3}\theta_{8} + \tilde{y}\theta_{5})^{N} 
    \\
    &= \gamma(x + x^2 + x^{3}\theta_{8})^{N} + \text{(coeff)}\tilde{y}\theta_{5}(x+ x^{2} + x^{3}\theta_{8})^{N-1} + \tilde{y}^{2}\theta_{9} .
\end{align*}


\noindent
Adding these equations and recalling that $\tilde{Y} = Y + \gamma {X}^{N}$, we find that 

\begin{align*}
    \tilde{Y} &= -\tilde{y}(1-\lambda x + \tilde{y}\theta_{1} + x^{2}\theta_{6}) + (2\gamma+\theta_{3})x^{N} + (\theta_{7} + \theta_{11})x^{N+1} + \theta_{10}x^{N-1}\tilde{y} + \theta_{9}\tilde{y}^{2}, 
\end{align*} 
i.e.,
    \begin{align*}
        \tilde{Y} &= -\tilde{y}(1-\lambda x + \tilde{y}[\theta_{1}-\theta_{9}] + x^{2}[\theta_{6}-x^{N-3}\theta_{10}]) + (2\gamma+\theta_{3})x^{N} + (\theta_{7} + \theta_{11})x^{N+1} 
        \\ 
    &= -\tilde{y}(1 - \lambda x + \tilde{y}\theta_{12} + x^{2}\theta_{13}) + (2\gamma + \theta_{3})x^{N} + \theta_{14}x^{N+1}. 
\end{align*}
(Recall $N \ge 3$.)

Now, $\theta_{3} = \beta + x\theta_{15} + \tilde{y}\theta_{16}$ for some number $\beta$ and some smooth functions $\theta_{15}, \theta_{16}$.  So, 
\begin{align*}
    \tilde{Y} &= -\tilde{y}(1 - \lambda x + \tilde{y}\theta_{12} + x^{2}\theta_{13} - x^{N}\theta_{16}) + (2\gamma + \beta)x^{N} + (\theta_{15} + \theta_{14})x^{N+1} .
    \\
    \mbox{Hence,}\    \tilde{Y} 
    &= -\tilde{y}(1 - \lambda x + \tilde{y}\theta_{12} + x^{2}\theta_{17}) + (2\gamma + \beta)x^{N} + \theta_{18}x^{N+1} .
\end{align*}


\noindent
Picking $\gamma = -\frac{\beta}{2} $, we kill the $(2\gamma + \beta)x^{N}$ term, leaving us with:

\begin{align} 
\label{eqn:Yt}
    \tilde{Y} = -\tilde{y}(1 - \lambda x + \tilde{y}\theta_{12} + x^{2}\theta_{17}) + x^{N+1}\theta_{18} .
\end{align}
Also,
\begin{align*}
    X = x + x^{2} + x^{3}\theta_{4} + [\tilde{y} - \gamma x^{N}]\theta_{5} ,
\end{align*}
so, 
\begin{align*}
    X = x + x^{2} + x^{3}[\theta_{4} - \gamma x^{N-3}\theta_{5}] + \tilde{y}\theta_{5} 
\end{align*}
i.e., 
\begin{align} 
\label{eqn:Xc} 
    X = x + x^{2} + x^{3}\theta_{19} + \tilde{y}\theta_{5} .
\end{align}

\noindent
(Again, recall that $N \ge 3$.)

Equations \eqref{eqn:Yt} and \eqref{eqn:Xc}  show that the map $\tilde{\Phi}: (x, \tilde{y}) \longmapsto (X, \tilde{Y})$ has the same form as $\Phi: (x, y) \longmapsto (X, Y)$, but 
with $N$ replaced by $(N+1)$.

\subsection{Conditions Preserved by Our Map}
Thanks to the preceding section, we may suppose our map has the form $\Phi : (x, y) \mapsto (X, Y) $, where 
\begin{equation} 
 \begin{aligned}
    X &= x + x^2 + x^{3}\theta_{A}(x, y) + y\theta_{B}(x,y) ,
    \\[1ex] 
    Y &= -y(1 - \lambda x + y\theta_{C}(x, y) + x^{2}\theta_{D}(x, y)) + x^{N+100}\theta_{E}(x, y) 
\end{aligned}
 \label{equ:IA}
\end{equation} 
with $\theta$'s smooth, $\lambda > 0$, N large. We fix $N \geq 100$ until further notice. We are only interested in $\Phi(x,y)$ when $x \geq 0$.

Let $K_{1}, ... , K_{N-10}$ be constants to be picked later. 
Initially, C, C', c, etc. will denote constants independent of the K's. Later, for each $m$, there will come a time when we will have picked $K_{1}, K_{2}, \ldots, K_{m}$, but have not yet picked $K_{m+1}, \ldots, K_{N-10}$. At that point, C, C', c, etc. will denote constants that may depend on $K_{1}, ..., K_{m}$, but not on $K_{m+1}, K_{m+2} ...$.

We spell out our conditions for the constants $K_1,...,K_{N-10}, \delta$.

\begin{itemize}
  \item $K_1$ is greater than a large enough constant determined by $\Phi$.
  \item Each $K_{m} \ (m \geq 2)$ exceeds a large enough constant determined by $\Phi, K_1,..., K_{m-1}$.
  \item $\delta$ is less than a small enough constant determined by $\Phi, K_1,..., K_{\bar{N}-10}$. 
\end{itemize}

\bigskip 

\begin{assertion}\label{assertion:1}
Suppose $|y| \leq x^{N}$. Let $(X,Y)=\Phi(x,y)$. Then $|Y|\leq X^N$, provided $x \in [0, 10\delta]$.
\end{assertion}

\begin{proof} 
Letting $C$ denote a constant determined by $\Phi$, we have 
\begin{align*}
    |Y| &\leq x^{N}(1 - \lambda x + Cx^{N} + Cx^{2}) + Cx^{N+100}
    \\
    &\leq x^{N} \ \ \mbox{for} \ x \in [0, 10\delta]
\end{align*}
while
\begin{align*}
    X &\geq x + x^{2} - Cx^{3} - Cx^{N} \geq x \text{, so $x^{N} \leq {X}^{N}$ for $x \in [0, 10\delta]$.}
\end{align*}
Therefore, $|Y| \leq x^{N} \leq X^{N} $ completing the proof of Assertion \ref{assertion:1}.
\end{proof}


\begin{assertion}\label{assertion:2}
Let $y=f(x)$ be a smooth function on $[0, x_{MAX}]$ with $0 < x_{MAX} \leq 10\delta$, and with
\begin{align*}
    |f(x)| &\leq x^{N} \ \ \mbox{for} \ x \in [0, x_{MAX}] ,
    \\
    |f'(x)| &\leq K_{1}x^{N-1} \ \ \mbox{for} \ x \in [0, x_{MAX}].
\end{align*}
Define $(X, Y) = \Phi(x, f(x))$ for $x \in [0, x_{MAX}]$. 
Then $\frac{dX}{dx} \geq 1$ for $x \in [0, x_{MAX}]$.
\end{assertion}

\begin{proof}
In fact,
\begin{align*}
    X = x + x^2 + x^{3}\theta_{A}(x, f(x)) + f(x)\theta_{B}(x, f(x)) .
\end{align*}
So,
\begin{align*}
    \frac{dX}{dx} = 1 + 2x + 3x^{2}\theta_{A}(x, f(x)) + x^{3}\theta_{A, x}(x, f(x)) + x^{3}\theta_{A, y}(x, f(x))f'(x) \\
    + f'(x)\theta_{B}(x, f(x)) + f(x)\theta_{B, x}(x, f(x)) + f(x)\theta_{B, y}(x, f(x))f'(x).
\end{align*}

\noindent
Because $|f(x)| \leq x^{N} $ and $|f'(x)| \leq K_{1}x^{N-1}$, all terms on the right-hand side  other than 1 and 2x are dominated by $Cx^{2}$ for $x \in [0, x_{MAX}]$ with C independent of $K_1$, since $x_{MAX} \leq 10\delta$.  (Recall, $\delta$ may depend on $K_1$.) So,
\[
\frac{dX}{dx} \geq 1 + 2x - Cx^2 \geq 1 ,
\] 
proving Assertion \ref{assertion:2}.

\end{proof}

\bigskip \bigskip


\begin{assertion}\label{assertion:3}
Let $f(x)$ be a smooth function on $[0, x_{MAX}]$ with $x_{MAX} \leq \delta $, and suppose that on that interval we have 
\begin{align*}
    |f(x)| \leq x^{N} \ \mbox{and}\  |f'(x)| \leq K_{1}x^{N-1}.
\end{align*} 
Then, if we pick $K_{1}$ large enough, it follows that
\begin{align*}
    \Phi(\{(x,y): x \in [0, x_{MAX}],  y=f(x)\})
   \  =  \ \{ (X, Y) : X \in [0, {X}_{MAX}], {Y} = F({X}) \}
\end{align*}
where $F$ is smooth and 
\begin{align*}
    |F({X})| \leq {X}^{N} ,  |F'({X})| \leq K_{1}{X}^{N-1} \text{ for } {X} \in [0, {X}_{MAX}].
\end{align*}
\end{assertion}

\begin{proof} 
Recall that until further notice, constants $C$ are determined by $\Phi$, independently of $K_{1},\ldots,K_{N-10},\delta$. By Assertion \ref{assertion:2}, we know that $\Phi(\{ (x, f(x)) : x \in [0, x_{MAX}] \})$ is the graph of a smooth function F on an interval $[0, {X}_{MAX}]$. By Assertion \ref{assertion:1}, we have $|F(X)|\leq X^N$ on $[0,X_{MAX}]$.

It remains only to estimate $|F'(X)|$. To do so, we note that, 
\begin{align*}
    F(X) &= {Y}, \text{with}
    \\
    {X} &= x + x^2 + x^3\theta_{A}(x, f(x)) + f(x)\theta_{B}(x, f(x))
    \\
    &\text{and}
    \\ 
    {Y} &= -f(x)(1-\lambda x + f(x) \theta_{C}(x, f(x)) + x^{2}\theta_{D}(x, f(x))) + x^{N+100}\theta_{E}(x, f(x)).
\end{align*}

\noindent
From the proof of Assertion \ref{assertion:2} and elementary calculus, we have 
\begin{enumerate}[a.]
\item
$\frac{d{X}}{dx} \geq 1 + 2x - Cx^{2} \geq 1$; $X \geq x$,
\item 
$F'(X)\frac{d{X}}{dx}=\frac{d{Y}}{dx}$, and 
\\[.4ex]
\item 
\vspace{-.4cm} 
\begin{align*} 
\frac{d{Y}}{dx} & = -f'(x)\{1 - \lambda x + f(x)\theta_{C}(x, f(x)) + x^{2}\theta_{D}(x, f(x))\}
\\
& - f(x)\{ -\lambda + f'(x)\theta_{C}(x, f(x)) + f(x)\theta_{C,x}(x, f(x)) + f(x)\theta_{C, y}(x, f(x))f'(x) 
\\
& + 2x\theta_{D}(x, f(x)) + x^{2}\theta_{D, x}(x, f(x)) + x^{2}\theta_{D, y}(x, f(x))f'(x) \} 
\\
& + (N+100)x^{N+99}\theta_{E}(x, f(x)) + x^{N+100}\theta_{E, x}(x, f(x)) + x^{N+100}\theta_{E, y}(x, f(x))f'(x) .
\end{align*} 

\end{enumerate}

\bigskip 

\noindent
Because $|f(x)| \leq x^{N}$,  $|f'(x)| \leq K_{1}x^{N-1}$, $x \leq x_{MAX} \leq \delta $ (and $\delta$ may depend on $K_{1}$), all the terms inside the second pair of curly brackets in (c) are dominated by C, and therefore

\begin{align*}
    |f(x) \{ -\lambda + f'(x)\theta_{C}(x, f(x)) + f(x)\theta_{C,x}(x, f(x)) + f(x)\theta_{C, y}(x, f(x))f'(x) 
\\
 + 2x\theta_{D}(x, f(x)) + x^{2}\theta_{D, x}(x, f(x)) + x^{2}\theta_{D, y}(x, f(x))f'(x) \}| \ \leq \ Cx^{N}.
\end{align*}

\noindent
Also because,
\begin{align*}
    |f(x)| \leq x^{N}, |f'(x)| \leq K_{1}x^{N-1} \text{ and } 0 \leq x \leq x_{MAX} \leq \delta,
\end{align*}

\noindent
we have 
\begin{align*}
    |f(x)\theta_{C}(x, f(x)) + x^{2}\theta_{D}(x, f(x))| \leq Cx^{2}
\end{align*}
and 
\begin{align*}
    |(N+100)x^{N+99}\theta_{E}(x, f(x)) + x^{N+100}\theta_{E,x}(x, f(x)) + x^{N+100}\theta_{E, y}(x, f(x))f'(x)| \leq Cx^{N+99}.
\end{align*}

\noindent
Therefore, (c) gives 
\begin{align*}
   \left| \frac{d{Y}}{dx} \right| &\leq |f'(x)|[1 - \lambda x + Cx^{2}] + Cx^{N}
    \\
    &\leq K_{1}x^{N-1}[1 - \lambda x + Cx^{2}] + Cx^{N} 
    \\
    &= K_{1}x^{N-1}[1-(\lambda - \frac{C}{K_{1}})x + Cx^{2}].
\end{align*}


\noindent
We pick $K_{1}$ large enough that $\lambda > \frac{C}{K_{1}}$. From now on, $K_{1}$ is fixed and constants C may depend on $K_{1}$. Then we have,

\begin{align*}
    \left |\frac{d{Y}}{dx} \right| \leq K_{1}x^{N-1} \text{ since } x \leq \delta \text{ and } \delta \text{ may depend on } K_{1}.
\end{align*}
Together with (a) and (b), this implies that 
\begin{align*}
    |F'({X})| \leq K_{1}x^{N-1} \leq K_{1}X^{N-1},
\end{align*}
completing the proof of Assertion \ref{assertion:3}. 
\end{proof}

\bigskip 
\noindent
We note the relationship of $x_{MAX}$ to ${X}_{MAX}$. We have 
\begin{align*}
    {X}_{MAX} = x_{MAX} + x_{MAX}^{2} + x^{3}_{MAX}\theta_A(x_{MAX}, f(x_{MAX})) + f(x_{MAX})\theta_{B}(x_{MAX}, f(x_{MAX})) 
\end{align*}
with $|f(x_{MAX})| \leq x^{N}_{MAX}$, so
\begin{align*}
    |{X}_{MAX} - (x_{MAX} + x_{MAX}^{2})|  \leq  Cx^{3}_{MAX}. 
\end{align*}


\bigskip 

\begin{assertion}\label{assertion:4}
Fix $\bar{m}\geq 2$, $\bar{m}\leq N-10$. Suppose the smooth function $f$ satisfies $|f(x)|\leq x^N$  and $|f^{(m)}(x)| \leq K_m x^{N-m}$ on $[0,x_{MAX}]$ for $m=1,...,\bar{m}$ with $x_{MAX} \leq \delta$. Define $F({X})$ on $[0, {X}_{MAX}]$ as in Assertion \ref{assertion:3}, and suppose 
\begin{itemize} 
\item $|F^{(m)}({X})| \leq K_{m} {X}^{N-m}$ for $m=1,...,\bar{m}-1$,  and 
\item $|F({X})|\leq {X}^N$
\end{itemize} 
on $[0, {X}_{MAX}]$.
\bigskip
\noindent
Here, $K_1,...K_{\bar{m}-1}$ have already been picked, but we have not yet picked $K_{\bar{m}}$. If $K_{\bar{m}}$ is large enough, then the above hypotheses imply the estimate, $|F^{(\bar{m})}({X})| \leq K_{\bar{m}} {X}^{N-\bar{m}}$.
\end{assertion}
\bigskip 

\begin{proof}
In this proof, $C, c, C'$ etc. denote constants determined by $\Phi, K_1,..., K_{\bar{m}-1}$. To prove Assertion \ref{assertion:4}, we differentiate the equation $F({X})={Y}$ \ \  $\bar{m}$ times with respect to $x$, where 
\begin{align*}
{X} &= x + x^2 +x^3 \theta_{A}(x, f(x)) + f(x) \theta_{B}(x,f(x)) , 
\\
{Y} &=-f(x) [ 1- \lambda x + f(x) \theta_{C}(x, f(x)) + x^2\theta_{D}(x,f(x))] + x^{N+100}\theta_{E}(x,f(x)).
\end{align*}
\noindent
Note that for $\theta(x,y)$ smooth, and for $1 \leq p \leq \bar{m}$, the quantity $(\frac{d}{dx})^{p} \theta(x, f(x))$ is a sum of terms $(\partial_{x}^{\alpha}\partial_{y}^{\beta} \theta) |_{(x, f(x))}\cdot\left[\prod\limits_{\nu=1}^{\beta}(\frac{d}{dx})^{p_{\nu}}f(x)\right]$ with each $p_{\nu} \geq 1$ and $\alpha + \sum\limits_{\nu}p_{\nu}=p$. 
The above term  is bounded by $C\cdot\prod\limits_{\nu=1}^{\beta}(K_{p_{\nu}}x^{N-p_{\nu}}) \leq C$, where the last estimate holds because $x \leq x_{MAX} \leq \delta$, and $p_{\nu} \leq \bar{m} \leq N-10$.  (Recall that $\delta$ is assumed to be less than a small constant determined by the $K_{p}$.)
\\
Therefore, 
\begin{align*}
    \left| \left(\frac{d}{dx}\right)^{p}\theta(x, f(x)) \right| \leq C \text{ for } 1 \leq p \leq \bar{m}. 
\end{align*}

\bigskip

\noindent
Together with the estimates we assumed for $\left(\frac{d}{dx}\right)^{p}f(x) \text{ where }(0 \leq p \leq \bar{m})$, this yields the following results  :

\begin{align*}
\left|\left(\frac{d}{dx}\right)^{p}{X}\right| \leq C \text{ for } 1 \leq p \leq \bar{m} ,
\end{align*}
\begin{align*}
\left|\left(\frac{d}{dx}\right)^{p}[(f(x))^2\theta_{C}(x, f(x))]\right| \leq (K_{\bar{m}} + C) \cdot C x^{2N-p}  \text{ for } 1 \leq p \leq \bar{m} , 
\end{align*}
\begin{align*}
\left|\left(\frac{d}{dx}\right)^{p}[x^{N+100}\theta_{E}(x, f(x))]\right| \leq Cx^{N+100-p}  \text{ for } 1 \leq p \leq \bar{m} , 
\end{align*}
\begin{align*}
\left|\left(\frac{d}{dx}\right)^{p}[f(x)\cdot x^2\theta_{D}(x, f(x))]\right| \leq C \cdot (K_{\bar{m}} + C) x^{N+2-p}  \text{ for } 1 \leq p \leq \bar{m} . 
\end{align*}
\noindent
So,
\begin{align} 
   \left(\frac{d}{dx}\right)^{\bar{m}}{Y} &  = -f^{(\bar{m})}(x)\cdot [1-\lambda x] + \lambda \bar{m} f^{(\bar{m}-1)}(x) + \text{Error}_{1} 
    \label{eqn:Smiley1} 
    \\ 
 \mbox{where} \ & |\text{Error}_{1}| \leq x^{N+1-\bar{m}} .
 \nonumber 
\end{align}

\bigskip
\noindent
(Here, we use the fact that $x \leq x_{MAX} \leq \delta$, where $\delta$ is less than a small constant depending on the K's.) Next, note that 
\begin{align*}
    \left(\frac{d}{dx}\right)^{\bar{m}}F({X})=F^{(\bar{m})}({X})\cdot \left(\frac{d{X}}{dx}\right)^{\bar{m}} + \sum\limits_{\substack{p \leq \bar{m} -1\\ r_1 + ...+ r_p= \bar{m}\\ \mbox{\footnotesize each} \ r_{p} \geq 1}}\text{(coeffs)} \ F^{(p)}({X})\cdot \prod\limits_{\nu=1}^p \left[\left(\frac{d}{dx}\right)^{r_{\nu}}{X}\right]. 
\end{align*}

\noindent
By our assumptions on $F^{(p)}({X})$ for $p \leq \bar{m}-1$, together with our estimates for $[(\frac{d}{dx})^r {X}]$ when $r \leq \bar{m}$, we therefore have 
\begin{align}
    \left(\frac{d}{dx}\right)^{\bar{m}}F({X}) & = F^{(\bar{m})}({X})\cdot \left(\frac{d{X}}{dx}\right)^{\bar{m}} + \text{Error}_{2}
    \label{eqn:Smiley2} 
        \\ 
 \mbox{where} \ & |\text{Error}_{2}| \leq C X^{N-\bar{m}+1} .
 \nonumber 
\end{align}

\noindent
Because $(\frac{d}{dx})^{\bar{m}}F({X}) = (\frac{d}{dx})^{\bar{m}}{Y}$, it follows from \eqref{eqn:Smiley1} and \eqref{eqn:Smiley2} that 
\begin{align*}
F^{(\bar{m})}({X})\cdot \left(\frac{d{X}}{dx}\right)^{\bar{m}} & =-f^{(\bar{m})}(x)\cdot (1-\lambda x) + \lambda \bar{m} f^{(\bar{m}-1)}(x) + \text{Error}_{3} 
        \\ 
 \mbox{where} \ & |\text{Error}_{3}| \leq C X^{N-\bar{m}+1}. 
\end{align*}
(Recall that $0 \leq x \leq {X}$.)

\noindent
Because $|f^{(\bar{m}-1)}(x)| \leq C \cdot x^{N-\bar{m}+1}$ (recall that $K_{\bar{m}-1}$ is a constant C), it therefore follows that 
\begin{align*}
F^{(\bar{m})}({X})\cdot \left(\frac{d{X}}{dx}\right)^{\bar{m}} & =-f^{(\bar{m})}(x)\cdot (1-\lambda x) + \text{Error}_{4} 
        \\ 
 \mbox{where} \ & |\text{Error}_{4}| \leq C x^{N-\bar{m}+1} .
\end{align*}
\noindent
(Here, we use the estimate $X \leq 2x$.)

\bigskip 

\noindent
Consequently,
\begin{align*}
|F^{(\bar{m})}({X})|\cdot \left(\frac{d{X}}{dx}\right)^{\bar{m}} &\leq |f^{(\bar{m})}(x)|\cdot (1-\lambda x) + C x^{N-\bar{m}+1} 
\\
&\leq K_{\bar{m}} x^{N-\bar{m}}(1-\lambda x) + C x^{N-\bar{m}+1} 
\\
&= K_{\bar{m}} x^{N-\bar{m}} \left(1-\lambda x + \frac{C}{K_{\bar{m}}}x \right)
\\
&\leq K_{\bar{m}} x^{N-\bar{m}},
\end{align*}

\noindent
provided we pick $K_{\bar{m}} \geq \frac{C}{\lambda}$, which we now do. Thus,
\begin{align*}
    |F^{(\bar{m})}({X})|\cdot \left(\frac{d{X}}{dx}\right)^{\bar{m}} \leq K_{\bar{m}} x^{N-\bar{m}} \leq K_{\bar{m}} {X}^{N-\bar{m}}.
\end{align*}

\noindent
We have seen that $\frac{d{X}}{dx} \geq 1$, hence the above estimate implies that $|F^{(\bar{m})}({X})| \leq K_{\bar{m}} {X}^{N-\bar{m}}$, completing the proof of Assertion \ref{assertion:4}. 
\end{proof}

\bigskip


\begin{assertion}\label{assertion:5}
For suitable constants $K_1,..,K_{N-10}$, and small enough $\delta$, the following holds. Let $\Gamma = \{(x, f(x)) : x \in [0, x_{MAX}]\}$ with $f$ smooth and 
\begin{align*}
    0 &< x_{MAX} \leq \delta , 
    \\
    |f(x)| &\leq x^N , 
    \\
    |f^{(m)}(x)| &\leq K_{m} x^{N-m} \ \  \mbox{for}\ m=1,..., N-10 . 
\end{align*}

\noindent
Then $\Phi(\Gamma)=\{({X}, F({X})): {X} \in [0, {X}_{MAX}]\}$ for a smooth function F that satisfies $|F({X})| \leq {X}^N$ and $|F^{(m)}({X})| \leq K_{m} {X}^{N-m}$ for $m=1,..., N-10$. Moreover, $|{X}_{MAX} - (x_{MAX} + x^{2}_{MAX})| \leq C x^{3}_{MAX}$. 
\end{assertion}

\bigskip 

\begin{proof}
We have already seen that $\Phi(\Gamma)=\{({X}, F({X})): {X} \in [0, {X}_{MAX}]\}$ with ${X}_{MAX}$ satisfying the estimate in Assertion \ref{assertion:5}. It remains only to check the estimates asserted for F.

We have seen that $|F({X})| \leq {X}^N$ and that $|F'({X})| \leq K_{1}{X}^{N-1}$. The desired estimates for $F^{(m)}({X}) \ \ \ (2 \leq m \leq N-10)$ follow from Assertion \ref{assertion:4} by an obvious induction on m. This completes the proof of Assertion \ref{assertion:5}.
\end{proof}

\noindent


\subsection{Shadowing}
In this section, we prove Lemma \ref{lemma:1.2}, stated in the introduction. We keep our assumption that $\Phi$ has the form (\ref{equ:IA}) for a fixed $N \geq 100$. Constants denoted $C$, $c$, $C'$, etc. may now depend on $\Phi, K_{1}, ..., K_{N-10}$, but not on $\delta$. Recall that we set $K_{0}=1$ in the introduction.

Assume 
\begin{equation}
  0 < x \leq \delta, \; |y| \leq x^N.
    \label{eqn:-1}
\end{equation}
Suppose
\begin{equation}
  \frac{|x-\hat{x}|+x^{-3}|y-\hat{y}|}{x^8} \leq 1.
    \label{eqn:0}
\end{equation}

\noindent
Define $(X,Y)$ and $(\hat{X},\hat{Y})$ by setting
\begin{equation}
    X= x+ x^2 + x^3\theta_{A}(x,y) + y\theta_B(x,y) , 
    \label{eqn:1}
\end{equation}
\begin{equation}
    \hat{X}= \hat{x} + \hat{x}^2 + \hat{x}^3\theta_{A}(\hat{x},\hat{y}) + \hat{y}\theta_B(\hat{x},\hat{y}) , 
    \label{eqn:2}
\end{equation}
\begin{equation}
    Y= -y \; (1- \lambda x + y \theta_C(x,y) +x^2 \theta_D(x,y)) + x^{N+100}\theta_E(x,y) , 
    \label{eqn:3}
\end{equation}
\begin{equation}
    \hat{Y}= -\hat{y} \; (1- \lambda \hat{x} + \hat{y} \theta_C(\hat{x},\hat{y}) +\hat{x}^2 \theta_D(\hat{x},\hat{y})) + \hat{x}^{N+100}\theta_E(\hat{x},\hat{y}) , 
    \label{eqn:4}
\end{equation}
i.e. $(X, Y)= \Phi(x,y)$ and $(\hat{X}, \hat{Y})= \Phi(\hat{x},\hat{y})$.

\noindent
We must show that 
\begin{equation}
    \frac{|X-\hat{X}|+X^{-3}|Y-\hat{Y}|}{X^8} \leq \frac{| x - \hat{x}| + x^{-3}|y-\hat{y}|}{x^{8}} \leq 1
    \label{eqn:5}
\end{equation}

\bigskip \bigskip 
\noindent
We start the proof of \eqref{eqn:5}. First of all, equations \eqref{eqn:-1} and \eqref{eqn:0} give 
\begin{equation}
    |\hat{y}| \leq C\hat{x}^8 \text{ and } |x-\hat{x}| \leq x^8. 
    \label{eqn:6}
\end{equation}
We note that 
\begin{equation}
    |\{-\lambda x + y \theta_C(x,y) + x^2 \theta_D(x,y)\} - \{-\lambda \hat{x} + \hat{y} \theta_C(\hat{x},\hat{y}) + \hat{x}^2 \theta_D(\hat{x},\hat{y})\}| \leq C[|x-\hat{x}|+|y-\hat{y}|]
    \label{eqn:7}
\end{equation}
and
\begin{eqnarray} 
\lefteqn{ |x^{N+100} \theta_E(x,y) -\hat{x}^{N+100}\theta_E(\hat{x},\hat{y})|
} \nonumber 
\\ 
&\leq&  |\hat{x}^{N+100}-x^{N+100}|\cdot |\theta_E(\hat{x},\hat{y})| +x^{N+100} |\theta_E(x,y) - \theta_E(\hat{x},\hat{y})| \nonumber   \\
    &\leq& x^N[|x-\hat{x}|+|y-\hat{y}|] . 
    \label{eqn:8}
\end{eqnarray} 

\noindent
Also,
\begin{equation}
    |(x+x^2)-(\hat{x} +\hat{x}^2)|=|x-\hat{x}|\cdot |1+x+\hat{x}| \leq [1 +2x +x^8] \cdot |x-\hat{x}|
    \label{eqn:8.5}
\end{equation}
by \eqref{eqn:6}.

\bigskip
\noindent
Similarly to \eqref{eqn:7} and \eqref{eqn:8}, we have
\begin{align} 
    |x^3 \theta_A(x,y) -\hat{x}^3\theta_A(\hat{x},\hat{y})| &\leq |x^3-\hat{x}^3|\cdot |\theta_A(\hat{x},\hat{y})|+ x^3 |\theta_A(x,y)-\theta_A(\hat{x},\hat{y})| \nonumber \\
    &\leq C x^2 |x-\hat{x}| + C x^3 [|x-\hat{x}|+|y-\hat{y}|] \nonumber \\
    &\leq C'x^2 |x-\hat{x}| + C' x^3 |y-\hat{y}| 
    \label{eqn:9}
\end{align}
and (thanks to \eqref{eqn:-1})
\begin{align} 
    |y \theta_B(x,y) -\hat{y}\theta_B(\hat{x},\hat{y})| &\leq |y-\hat{y}|\cdot |\theta_B(\hat{x},\hat{y})|+ |y|\cdot |\theta_B(x,y)-\theta_B(\hat{x},\hat{y})| \nonumber 
    \\
    &\leq C |y-\hat{y}| + C x^N [|x-\hat{x}|+|y-\hat{y}|] \nonumber 
\\
    &\leq Cx^N |x-\hat{x}| + C |y-\hat{y}|. 
    \label{eqn:10}
\end{align}
\noindent
We apply the above to estimate $|X-\hat{X}|$ and $|Y-\hat{Y}|$.

\bigskip 


\noindent
We have from \eqref{eqn:1}  and \eqref{eqn:2}  that 

\begin{align} 
    |X-\hat{X}| & \leq  |(x+x^2) - (\hat{x} + \hat{x}^2)| + |x^3 \theta_A(x,y)-\hat{x}^3\theta_A(\hat{x},\hat{y})| + |y \theta_B(x,y)-\hat{y}\theta_B(\hat{x},\hat{y})| \nonumber 
    \\ 
   & \leq [1 + 2x + x^8]|x-\hat{x}| + \{C'x^2|x - \hat{x}| + C'x^3|y - \hat{y}|\} + \{Cx^N|x-\hat{x}|+C|y - \hat{y}|\} 
   \nonumber 
   \\ 
   & \leq [1 + 2x + Cx^2] |x - \hat{x}| + C |y - \hat{y }| , 
    \label{eqn:11}
\end{align} 
where the second inequality follows by \eqref{eqn:8.5}, \eqref{eqn:9}, \eqref{eqn:10}.

From \eqref{eqn:3} and \eqref{eqn:4} we have 

\begin{align} 
|Y-\hat{Y}| \ \leq \ & |y - \hat{y}| \,  |\{1-\lambda \hat{x} + \hat{y} \theta_C (\hat{x}, \hat{y}) + \hat{x}^2 \theta_D(\hat{x}, \hat{y})\}| 
\nonumber 
\\ 
& + |y| \, |\{ - \lambda x  + y \theta_C(x,y) + x^2 \theta_D(x,y)\} - \{ -\lambda \hat{x} + \hat{y} \theta_c(\hat{x}, \hat{y}) + \hat{x}^2 \theta_D(\hat{x}, \hat{y}) \} |  
\nonumber 
\\ & + |x^{N+100} \theta_E(x,y) - \hat{x}^{N+100} \theta_E(\hat{x}, \hat{y})| 
\nonumber 
\\ 
\ \leq \ & |y - \hat{y}| + C x^N [|x-\hat{x}| + |y-\hat{y}|]  
\nonumber 
\end{align}

\noindent 
where we have used the fact that $\lambda > 0$ and 
\begin{align} 
|\hat{y} \theta_C(\hat{x}, \hat{y}) + \hat{x}^2 \theta_D (\hat{x}, \hat{y})| \leq C |\hat{y}| + C \hat{x}^2 \leq C'\hat{x}^2 
\nonumber 
\end{align}
 by \eqref{eqn:6}. 

\bigskip

\noindent
Consequently, 
 \begin{align} 
 |Y-\hat{Y}| \leq (1+C x^N) |y - \hat{y}| + C x^N |x - \hat{x}|.
 \label{eqn:14}
 \end{align}
Note also that 
\begin{align} 
 x \leq x + x^2 - C x^3 \leq X \leq x + x^2 + C x^3 ,
 \label{eqn:15}
\end{align}
thanks to \eqref{eqn:-1} and \eqref{eqn:1}. 
In particular $ X^{-3} \leq x^{-3}$, so \eqref{eqn:14} yields 
\begin{align*} 
     X^{-3} |Y - \hat{Y}| \leq C x^{N-3} |x-\hat{x}| + (1 + C x^N) x^{-3} |y - \hat{y}|.  
\end{align*}
Adding this to \eqref{eqn:11}, we find that 
\begin{align} 
    |X-\hat{X}| + X^{-3}|Y-\hat{Y}| & \leq [1 + 2x + Cx^2]|x-\hat{x}| + (1+Cx^N + C x^3)x^{-3}|y-\hat{y}|  \nonumber 
    \\ 
 & \leq [1+2x+Cx^2][|x-\hat{x}|+x^{-3}|y-\hat{y}|].
    \label{eqn:17}
\end{align}

From \eqref{eqn:15} we have also 
\begin{align} 
X^{-8} \leq (x + x^2 -Cx^3)^{-8} = x^{-8}(1+x-Cx^2)^{-8} \leq x^{-8}(1-8x + C'x^2).  
\label{eqn:18}
\end{align}

\bigskip 

\noindent
Multiplying \eqref{eqn:17} by \eqref{eqn:18} we have 

\begin{align*} 
\frac{|X-\hat{X}|+ X^{-3}|Y-\hat{Y}|}{X^8} & \leq (1-8x+C'x^2)(1+2x+Cx^2) \left( \frac{|x-\hat{x}|+x^{-3}|y-\hat{y}|}{x^8} \right) 
\\ 
&\leq \frac{|x-\hat{x}|+x^{-3}|y-\hat{y}|}{x^8}. 
\end{align*} 
In particular, recalling our assumption \eqref{eqn:0}, we see that \eqref{eqn:5} holds. This completes the proof of our shadowing result \eqref{eqn:5}, thus establishing Lemma \ref{lemma:1.2}. $\square$
\vspace{5mm}

\vspace{3mm}

\noindent 

\vspace{5mm}
\subsection{Approximately Invariant Curves}

\vspace{5mm}

Fix constants $K_1, ... , K_{N-10}$, $\delta$ as before. As in the introduction, we set $K_{0}=1$. Constants $C$ will depend on those K's, but not on $\delta$. Let $0 < \rho < \delta$ be a small number. Later, we will fix $\delta$ small enough, and let $ \rho \rightarrow 0^{+}$. We recall that $\delta$ is less than a small enough constant determined by $ \Phi$, $K_1$, ..., $K_{N-10}$.

Let $ x_{MAX}^{0} = \rho  $   and let $F_0(x) = 0$ on $[0, x_{MAX}^0 ]$. By induction on $\nu$, define $x_{MAX}^{\nu}$ and $F_{\nu}(x)$ on $[0, x_{MAX}^{\nu}]$ by setting $\Phi(\{(x, F_{\nu-1}(x)): x \in [0, x_{MAX}^{\nu-1}]\}) = \{(x, F_{\nu}(x)): x \in [0, x_{MAX}^{\nu}] \}$. We terminate the construction of $F_{\nu}$, $x_{MAX}^{\nu}$ as soon as we can no longer apply our Assertions \ref{assertion:1}--\ref{assertion:5} to keep going. That is, we pass from $x_{MAX}^{\nu-1}$, $F_{\nu-1}$ to $x_{MAX}^{\nu}$, $F_{\nu}$ provided $x_{MAX}^{\nu-1} \leq \delta$. If $x_{MAX}^{
\nu-1} > \delta $, then we stop.

\noindent
As long as the $x_{MAX}^{\nu}$, $F_{\nu}$ are well-defined, we have 
\[|F_{\nu}(x)| \leq x^N \text{ for } x \in [0,x_{MAX}^{\nu}]\]
and
\[|F_{\nu}^{(m)}(x)| \leq K_m x^{N-m} \text{ for } x \in [0, x_{MAX}^{\nu}], \ \ 1 \leq m \leq N - 10 . \] 

\noindent
Indeed, that holds for $\nu = 0$ since $F_0 \equiv 0$; and it then follows by induction thanks to Assertion \ref{assertion:5}. Note that $x_{MAX}^{\nu} \geq x_{MAX}^{\nu-1} + (x_{MAX}^{\nu-1})^2 - C(x_{MAX}^{\nu-1})^3$. As long as $x_{MAX}^{\nu-1} \leq \delta$ and $\delta$ is less than a small enough constant, we have $x_{MAX}^{\nu} \geq x_{MAX}^{\nu-1} + \frac{1}{2} (x_{MAX}^{\nu-1})^2$. Consequently, our induction on $\nu$ will eventually terminate, i.e, $x_{MAX}^{\nu} > \delta$ for some $\nu$. Let $\bar{\nu}$ denote the first $\nu$ for which $x_{MAX}^{\nu} > \delta$. Thus, our induction defines $F_0, F_1, ..., F_{\bar{\nu}}$  but then terminates. We have $\delta < x^{\bar{\nu}}_{MAX} \leq x^{\bar{\nu}-1}_{MAX} + (x^{\bar{\nu}-1}_{MAX})^{2} + C(x_{MAX}^{\bar{\nu}-1})^{3} \leq \delta + \delta^{2} + C\delta^{3} \leq 2\delta $.

Now suppose $\bar{x} \in [0, \frac{\delta}{2}]$ is given.
 Then $(\bar{x}, F_{\bar{\nu}}(\bar{x})) = \Phi^{\bar{\nu}}(\bar{x}_0, 0)$ for some $\bar{x}_{0} \in [0, \rho]$. (That's because, for any $\nu \in \{ 0, 1, ...., \bar{\nu}\}$, we have 
 $\Phi^{\nu} \{ (x,0): x \in [0, \rho ]\} = \{ (x, F_{\nu}(x)) : x \in [0,x_{MAX}^{\nu} ]\}$.) 
 
 \noindent
 There exists $\Tilde{x}_0 \in [0, \bar{x}_0]$ such that 
 \[ \Tilde{x}_0 + \Tilde{x}_0^2 + \Tilde{x}_0^3 \theta_A (\Tilde{x}_0,0) = \bar{x}_0.\]
 \noindent
 Then $\Phi(\Tilde{x}_0,0) = (\bar{x}_0, \Tilde{y}_0)$ with $|\Tilde{y}_0| \leq \Tilde{x}_0^N \leq \bar{x}_0^N \leq \rho^{N}$. (The first inequality here is immediate from (\ref{equ:IA})). Let $z_{\nu} = (x_{\nu},y_{\nu}) = \Phi^{\nu}(\bar{x}_0,0)$
 and 
 $ \hat{z}_{\nu} = (\hat{x}_{\nu},\hat{y}_{\nu}) = \Phi^{\nu}(\bar{x}_0,\Tilde{y}_0)$ 
 for $\nu = 0, ..., \bar{\nu}$. Let's estimate how close $\hat{z}_{\nu}$ is to $z_{\nu}$.

\noindent
 Note that 
\[ \frac{|x_0 - \hat{x}_0| + x_0^{-3}|y_0-\hat{y}_0|}{x_0^{8}} = \frac{|\tilde{y}_0|}{\bar{x}_0^{11}} \leq \bar{x}_0^{N-11} \leq \rho^{N-11}.\] 

\noindent
Repeatedly applying  Lemma~\ref{lemma:1.2}, we see that 

\[ \frac{|x_{\bar{\nu}}-\hat{x}_{\bar{\nu}}| + x_{\bar{\nu}}^{-3} |y_{\bar{\nu}}-\hat{y}_{\bar{\nu}}|  }{x_{\bar{\nu}}^8} \leq \rho^{N-11}.\]

\noindent
In particular, 
\[ |x_{\bar{\nu}}-\hat{x}_{\bar{\nu}}|, |y_{\bar{\nu}}-\hat{y}_{\bar{\nu}}| \leq \rho^{N-11}, \ \mbox{so} \]

\[|z_{\bar{\nu}}-\hat{z}_{\bar{\nu}}| \leq C \rho^{N-11}. \]

\noindent
That is,
\[|\Phi^{\bar{\nu}}(\bar{x}_0,0) - \Phi^{\bar{\nu}}(\bar{x}_0, \Tilde{y}_0)| \leq C \rho^{N-11}.\]
\newline

\noindent
Recall that
\begin{align*}
    \Phi^{\bar{\nu}}(\bar{x}_{0}, 0) = (\bar{x}, F_{\bar{\nu}}(\bar{x}))
\end{align*}
and
\begin{align*}
    (\bar{x}_{0}, \tilde{y}_{0}) = \Phi(\tilde{x}_{0}, 0), \ \mbox{so}
\end{align*} 
\[ 
\Phi^{\bar{\nu}}(\bar{x}_{0}, \tilde{y}_{0}) 
= \Phi^{\bar{\nu}+1}(\tilde{x}_{0}, 0) 
= \Phi(\Phi^{\bar{\nu}}(\tilde{x}_{0}, 0)). 
\] 
Consequently,
\begin{align*}
    |(\bar{x}, F_{\bar{\nu}}(\bar{x})) - \Phi(\Phi^{\bar{\nu}}(\tilde{x}_{0}, 0))| \leq C \rho^{N - 11}. 
\end{align*}

\noindent
Now $\tilde{x}_{0} \in [0, \bar{x}_{0}] \subset [0, \rho]$, and %
\begin{align*}
    \Phi^{\bar{\nu}}(\{(x, 0): x \in [0, \rho]\}) &= \{(x, F_{\bar{\nu}}(x)) : x \in [0, x_{MAX}^{\bar{\nu}}]\},
\end{align*}
hence 
\begin{align*}
    \Phi^{\bar{\nu}}(\tilde{x}_{0}, 0) &= (\hat{x}, F_{\bar{\nu}}(\hat{x})) 
\end{align*}
for some  $\hat{x} \in [0, x^{\bar{\nu}}_{MAX}] \subset [0, 2\delta]$.
\noindent
Thus,
\begin{align}
    \label{eq:label}
    |(\bar{x}, F_{\bar{\nu}}(\bar{x})) - \Phi(\hat{x}, F_{\bar{\nu}}(\hat{x}))| \leq C\rho^{N-11}.
\end{align}
\noindent
We have $\Phi(\hat{x}, F_{\bar{\nu}}(\hat{x})) = (\hat{X}, \hat{Y})$, with $|\hat{X}-(\hat{x}+\hat{x}^{2})| \leq C |\hat{x}|^{3} + C|F_{\bar{\nu}}(\hat{x})| \leq C'|\hat{x}|^{3} $ and $|\hat{x}| \leq 2\delta$. Hence $\hat{X}\geq \hat{x}$.

On the other hand, (\ref{eq:label}) gives $|\hat{X}-\bar{x}|\leq C\rho^{N-11}$. Therefore, 
\[ \hat{x} \leq \bar{x} + C\rho^{N-11} \leq \frac{\delta}{2} + C\rho^{N-11} \leq \delta, 
\] provided $C\rho^{N-11} < \frac{\delta}{2}$. We have established the following result. 

\begin{assertion}
\label{assertion:6}
Suppose $\rho$ is less than a small enough positive constant determined by $\Phi, K_{1},\ldots,K_{N-10}, \delta$. Then, given $\bar{x}\in [0, {\frac{\delta}{2}}]$, there exists $\hat{x} \in [0, \delta]$ such that \[ 
\hat{x} \leq \bar{x} + C\rho^{N-11}
\] 
and 
\[ 
|(\bar{x}, F_{\bar{\nu}}(\bar{x})) - \Phi(\hat{x}, F_{\bar{\nu}}(\hat{x}))| \leq C\rho^{N-11}. 
\] 
Moreover, $F_{\bar{\nu}}$ satisfies 
\[ 
\left| \left(\frac{d}{dx}\right)^{m}F_{\bar{\nu}}(x) \right| \leq K_{m}x^{N-m} \ \  \mbox{for} \ 0 \leq m \leq N-10,  \ x \in [0, \delta].
\]
\end{assertion}

\bigskip 

\subsection{Passing to the Limit}
As before, we fix $N\geq100$ and suppose our map $\Phi$ has the form (\ref{equ:IA}). We fix the constants $K_{1}, ..., K_{N-10}, \delta$ and consider a sequence $\rho_{1}, \rho_{2}, ...$ of positive numbers tending to zero. For each $\rho_{j}$, we apply Assertion \ref{assertion:6}. 

\noindent
Thus, we obtain a sequence of functions $F_{j} \in C^{N-10}([0, \delta])$, with the following properties.

\begin{align}
    \label{eq:38}
    |F_{j}(x)| \leq x^{N} \text{ for } x\in[0, \delta].
\end{align}
\begin{align}
    \label{eq:39}
    \left| \left( \frac{d}{dx} \right)^{m}F_{j}(x)\right| \leq K_{m}x^{N-m} \text{ for } x\in[0, \delta], 1 \leq m \leq N-10 . 
\end{align}
 
\noindent  
Given $\bar{x} \in \left[0, \frac{\delta}{2}\right]$  there exists $\hat{x}_{j}\in [0, \delta] \cap \left[0, \bar{x}+C\rho_{j}^{N-11}\right]$ such that     
\begin{align}
    \label{eq:40}
    |(\bar{x}, F_{j}(\bar{x})) - \Phi(\hat{x}_{j}, F_{j}(\hat{x}_{j}))| \leq C\rho_{j}^{N-11}.
\end{align}    

\noindent    
By Ascoli's Theorem, we may pass to a subsequence to achieve for some $F \in C^{N-11}$ $\left(\left[0,\delta\right]\right)$ that
\begin{align}
    \label{eq:41}
    F_{j} \rightarrow F \text{ in } C^{N-11} \text{ norm.}
\end{align}

\noindent
From (\ref{eq:38}) and (\ref{eq:39}), we have
\begin{align}
    |F(x)| \leq x^{N} \text{ for } x\in[0, \delta]
\end{align}
and
\begin{align}
    \left|\left(\frac{d}{dx}\right)^{m}F(x)\right| \leq K_{m}x^{N-m} \text{ for } x\in[0, \delta], 1\leq m \leq N-11.
\end{align}

\noindent
Now let $\bar{x}\in \left[0, \frac{\delta}{2}\right]$, and let $\hat{x}_{j}$ be as in (\ref{eq:40}). Passing to a subsequence $\hat{x}_{j_{i}}(i=1,2,3,...)$ depending on $\bar{x}$, we may achieve,     

\begin{align}
    \label{eq:44}
    \hat{x}_{j_{i}} \rightarrow \hat{x} \text{ as } i\rightarrow\infty
\end{align}
with, 
   
\begin{align}
    \hat{x}\in[0, \bar{x}] \subset \left[0, \frac{\delta}{2}\right] . 
\end{align}   
   
\noindent
Thanks to (\ref{eq:41}) and (\ref{eq:44}), we have \[ 
(\hat{x}_{j_{i}}, F_{j_{i}}(\hat{x}_{j_{i}})) \rightarrow (\hat{x}, F(\hat{x})) 
\ \mbox{as} \ i\rightarrow\infty,
\] 
hence,
\begin{align}
    \label{eq:46}
    \Phi(\hat{x}_{j_{i}}, F_{j_{i}}(\hat{x}_{j_{i}})) \rightarrow \Phi(\hat{x}, F(\hat{x})) \text{ as } i \rightarrow\infty.
\end{align}   
We now have 
\begin{eqnarray}  
|(\bar{x}, F(\bar{x})) - \Phi(\hat{x}, F(\hat{x}))| & \leq & 
|(\bar{x}, F(\bar{x})) - (\bar{x}, F_{j_{i}}(\bar{x}))| +
    |(\bar{x}, F_{j_{i}}(\bar{x})) - \Phi(\hat{x}_{j_{i}}, F_{j_{i}}(\hat{x}_{j_{i}}))| 
       \nonumber 
    \\
   & & + | \Phi(\hat{x}_{j_{i}}, F_{j_{i}}(\hat{x}_{j_{i}})) - \Phi(\hat{x}, F(\hat{x}))|.
    \label{eq:47} 
\end{eqnarray}

\noindent
The three terms on the right in (\ref{eq:47}) all tend to zero as $i\rightarrow\infty$, thanks to (\ref{eq:40}), (\ref{eq:41}), and (\ref{eq:46}).  Hence, $(\bar{x}, F(\bar{x}))=\Phi(\hat{x}, F(\hat{x}))$. 
\noindent
We have therefore proven the following. 

\begin{assertion}\label{A7}
Let $N\geq100$, and suppose $\Phi$ has the form (\ref{equ:IA}). Then there exist $\delta>0$ and $F\in C^{N-11}([0, \delta])$ with the following properties,
\begin{itemize}
    \item $|F(x)| \leq x^{N} \text{ for } x \in [0, \delta]$.
    \item $\left| \left( \frac{d}{dx}\right)^{m}F(x) \right| \leq K_{m}x^{N-m}$ for $x\in[0,\delta], 1\leq m\leq N-11$.
    \item Given $\bar{x}\in \left[0, \frac{\delta}{2}\right]$ there exists $\hat{x}\in [0, \bar{x}]$ such that $(\bar{x}, F(\bar{x}))=\Phi(\hat{x}, F(\hat{x}))$.
\end{itemize}
\end{assertion}

\bigskip 

We now pass from the setting of maps (\ref{equ:IA}) back to our original coordinates, in which our map $\Phi$ has the form (\ref{eqn:m2}). 

Recall that we pass from (\ref{eqn:m2}) to (\ref{equ:IA}) by repeatedly making coordinate changes of the form $(x, y) \rightarrow (x, y + \gamma x^{n}), (X, Y) \rightarrow (X, Y+\gamma X^{n})$ with $n\geq 3$. From Assertion \ref{A7}, we therefore read off the following conclusion.

\begin{assertion}\label{A8}
Let $\Phi$ be a mapping of the form (\ref{eqn:m2}), and let $N\geq 100$ be given. Then there exist positive constants $\delta_{N}, C_{N}$ and a function $F_{N}\in C^{N}([0, \delta_{N}])$ with the following properties. 

\begin{enumerate}[(I)]
\item  \textit{Tangency:} $|F_{N}(x)| \leq C_{N}x^{3}$ for $x\in [0, \delta_{N}]$ . 

\item \textit{Invariance:} Given $\bar{x}\in[0, \delta_{N}]$ there exists $\hat{x}\in [0, \bar{x}]$ such that $(\bar{x}, F_{N}(\bar{x})) = \Phi(\hat{x}, F_{N}(\hat{x}))$ . 

\end{enumerate}
\end{assertion}

\bigskip 

\noindent
So far, $F_{N}$ and $\delta_{N}$ may depend on N. In the next section, we remedy this defect.

\subsection{Uniqueness}
\noindent
We prove the following local uniqueness result.
\noindent
\begin{assertion}\label{A9}
Let $\Phi$ be as in (\ref{eqn:m2}), let $N\geq100$, and let $F_{N}\in C^{N}([0,\delta_{N}])$ be as in Assertion \ref{A8}. Suppose $\tilde{F}:[0, \tilde{\delta}] \rightarrow \mathbb{R}$ satisfies 
\begin{itemize} 
\item $ x^{-2/3}\tilde{F}(x) \rightarrow 0 \text{ as }x\rightarrow 0^{+}$ \\ 
and 
\item  For every  $\bar{x} \in [0, \tilde{\delta}]$  there exists  $\hat{x} \in [0, \bar{x}]$  such that  $(\bar{x}, \tilde{F}(\bar{x})) = \Phi(\hat{x}, \tilde{F}(\hat{x}))$ . 
\end{itemize}

\noindent
Then for some small positive $\delta \leq \min (\delta_{N}, \tilde{\delta})$ we have $\tilde{F} = F_{N} \text{ on } [0, \delta]$. 
\end{assertion}

\begin{proof}
By making a change of coordinates, 
\begin{align*}
    y^{\#} &= y - F_{N}(x), x^{\#} = x , 
    \\
    Y^{\#} &= Y - F_{N}(X), X^{\#} = X,
\end{align*}
we may assume without loss of generality that $F_{N} = 0 \text{ on } [0, \delta_{N}]$. (However, $\Phi$ is now merely $C^{N}$, not $C^{\infty}$.) We must show that $\tilde{F}(x) = 0$ for small positive $x$. Thanks to the invariance condition in Assertion \ref{A8}, with $F_{N}=0$, our map $\Phi$ has the form $(x, y) \mapsto (X, Y)$ with
\begin{align*}
    X &= x + x^{2} +\mu xy + O(|(x, y)|^{3})
    \\
    Y &= -y(1-\lambda x + O(|(x,y)|^{2})).
\end{align*}
\noindent
Hence, $\Phi^{-2}$ has the form $(x, y) \mapsto (X, Y)$ with 
\begin{equation}
    \begin{aligned}
     X &= x - 2x^{2} + O(|(x, y)|^{3}) 
     \\ 
     Y &= y(1+2\lambda x + O(|(x,y)|^{2})).
    \end{aligned}
    \label{eq:50}
\end{equation}

\noindent
Note that the term $\mu x y$ above contributes only $O(|(x,y)|^3)$ to $\Phi^{-2}(x,y)$. We study $\Phi^{-2}(x, \tilde{F}(x))$ for small positive $x$. 

Since $|\tilde{F}(x)|=o(x^{2/3})$, we have
\begin{align*}
    |(x, \tilde{F}(x))|^{3} & = o(x^{2}) 
    \end{align*}
and 
    \begin{align*}
    |(x,\tilde{F}(x))|^{2} & = o(x).
\end{align*}
\noindent
Consequently, (\ref{eq:50}) and the invariance property of $\tilde{F}$ together imply for $x>0$ small enough:

\begin{align}
    \label{eq:51}
    \Phi^{-2}(x, \tilde{F}(x)) = (\hat{x}, \tilde{F}(\hat{x})) 
\end{align}
\begin{align}
\intertext{with}  
  \label{eq:52}
    0 < \hat{x} < x- \frac{1}{2}x^{2} 
\end{align}
\begin{align}
\shortintertext{and} 
    \label{eq:53}
    |\tilde{F}(\hat{x})| \geq |\tilde{F}(x)|.
\end{align}

\noindent
Now suppose that for some small enough positive $x_{0}$, we have $\tilde{F}(x_{0})\neq0$. Repeatedly applying (\ref{eq:51}), (\ref{eq:52}), (\ref{eq:53}), we learn that $\Phi^{-2 \nu}(x_{0}, \tilde{F}(x_{0})) = (x_{\nu}, \tilde{F}(x_{\nu}))$, with $x_{\nu} \rightarrow 0$ as $\nu \rightarrow \infty$, but $|\tilde{F}(x_{\nu})| \geq |\tilde{F}(x_{0})| > 0$ for all $\nu$. This contradicts our hypothesis $|\tilde{F}(x)|=o(x^{2/3})$. Thus, $\tilde{F}(x_{0})=0$ for all small enough $x_{0}>0$, completing the proof of Assertion \ref{A9}. 

\end{proof}

\subsection{Endgame}
Let $\Phi$ be as in the statement of Theorem \ref{thm1.1}. For each $N \geq 100$, let $F_N \in C^N([0, \delta_N])$ be as in Assertion \ref{A8}.

Assertion \ref{A9} tells us that $F_N=F_{N'}$ in an interval $[0, \delta_{(N,N')}]$ for all $N, N' \geq 100$. In particular, $F_N=F_{100}$ on an interval $[0, \delta_{(N,100)}]$ for each $N \geq 100$. 

Consequently, $F_{100} \in C^N([0, \delta_{(N,100)}])$ for each such $N$. Repeatedly applying the invariance condition in Assertion \ref{A8} to $F_{100}$, we learn that, for any $\nu \geq 1$, the graph $\Gamma=\{(x,F_{100}(x)): x \in [0, \delta_{100}]\}$ is equal to the image of the graph $\{(x,F_{100}(x)): x \in [0, \hat{\delta}_{(\nu)}]\}$ under the map $\Phi^{\nu}$, for some $\hat{\delta}_{(\nu)}>0$. 

We have $\hat{\delta}_{(\nu +1)} \leq \hat{\delta}_{(\nu)} - \frac{1}{2}(\hat{\delta}_{(\nu)})^2$, hence $\hat{\delta}_{(\nu)} \rightarrow 0 $ as $\nu \rightarrow \infty $. Taking $\nu$ so large that $\hat{\delta}_{(\nu)} < \delta_{(N,100)}$, we see that $\Gamma$ is the image of a $C^N$ curve under the smooth map  $\Phi^{\nu}$. Therefore, $F_{100} \in C^N([0, \delta_{100}])$ for all $N \geq 100$. 

Thus, $F_{100} \in C^{\infty}([0, \delta_{100}])$. Together with Assertion \ref{A8} for $F_{100}$, this proves the existence claimed in Theorem \ref{thm1.1}. Finally, the uniqueness claimed in Theorem \ref{thm1.1} is precisely Assertion \ref{A9}. 

The proof of Theorem \ref{thm1.1} is complete. 
$\blacksquare$



\appendix
\section*{Appendix}

Theorem 1.1 follows from Theorem 2.1 in \cite{baldoma2007} applied to the inverse map of $\Phi^2$.
Indeed, Theorem 2.1 provides, under appropriate conditions, a stable manifold of the origin, tangent to the $x$ axis, for a map $\Psi$ such that $\Psi(0,0) =(0,0)$ and
$D\Psi(0,0) =\begin{pmatrix} 1 &0 \\0 &  1 \end{pmatrix}$.
Given $\Phi$ as in (1) in your preprint, we have
$$
\Phi^2(x,y)
=
\begin{pmatrix} x+2x^2 + O(|(x,y)|^3) \\y-2\lambda xy + yO(|(x,y)|^2) + O(|(x,y)|^4) \end{pmatrix}
$$
(An important point here is that the second component of $
\Phi^2(x,y)$ has no term of the form $cx^3$.)

Then 	
$$
\Psi(x,y)
=\begin{pmatrix} \Psi_1(x,y) \\ \Psi_2(x,y)  \end{pmatrix}
=\Phi^{-2}(x,y) =  \begin{pmatrix} x-2x^2 + O(|(x,y)|^3) \\ y+2\lambda xy + yO(|(x,y)|^2) + O(|(x,y)|^4) \end{pmatrix}.
$$
Taking $ F= \Psi $, $N=M=2$	in the statement of Theorem 2.1
we check that
$$
\frac{\partial^2 \Psi_1 }{\partial x^2} (0,0) = -2 <0, \qquad
\frac{\partial^2 \Psi_2 }{\partial x^2} (0,0) = 0, \qquad
\frac{\partial^2 \Psi_1 }{\partial x \partial y} (0,0) = 2 \lambda >0.
$$	
Then there exists a $C^\infty $ map $K:[0,t_0] \to \R^2$ and a polynomial (of degree 3) $R: \R\to \R$ such that
$$
\Psi \circ K = K \circ R.
$$
The image of $K$ is the stable manifold of $\Psi$.
The polynomial approximation of $K$ given in Section 3 of \cite{baldoma2007} provides, using that there is no term of the form $cx^3$ in $\Psi_2$,
$$
K(t) = \begin{pmatrix}
t+O(t^3) \\ O(t^3)
\end{pmatrix}, \qquad R(t) = t-2t^2 +d t^3, \quad d\in \R.
$$
$K$ is not unique but its image is the graph of a unique function $\varphi $ (see Remark 2.3).
Also $\varphi(x) = O(x^3)$.
The uniqueness is among all Lipschitz functions from $[0, x_0]$ to $\R$ satisfying $|\varphi(x)|\leq C |x|$ for arbitrary constant, changing if necessary the value of $x_0$. Our uniqueness statement is not exactly yours.  

Now we have that $\text{graph }  \varphi $ is invariant by $\Psi$. We take $\tilde \varphi$ such that
$$
\text{graph } \tilde \varphi =  \Phi^{-1}(\text{graph } \varphi)
$$
(in a slightly smaller domain). Since
$$
\Phi^{-1}  \begin{pmatrix}
	x \\ \varphi (x)
\end{pmatrix} =
 \begin{pmatrix}
x-x^2 + O(x^3) \\ -\varphi (x) -\lambda x \varphi (x) + O(x^3)
\end{pmatrix}
=
\begin{pmatrix}
	x-x^2 + O(x^3) \\ O(x^3)
\end{pmatrix}
$$
we have that $ \tilde \varphi (x)=O(x^3)$.

Moreover, since
$$
\Psi(\text{graph } \tilde \varphi ) = \Phi^{-3}(\text{graph } \varphi ) = \Phi^{-1} (\Psi(\text{graph }  \varphi ) )
\subset \Phi^{-1}(\text{graph }  \varphi ) = \text{graph } \tilde \varphi,
$$
$\text{graph } \tilde \varphi$ is also invariant by $\Psi$.
Then, by the uniqueness property, $\varphi=\tilde \varphi$
in the common domain.

\textbf{Final remark 1} In \cite{baldoma2020} we deal with invariant manifolds of arbitrary (finite) dimension. One could deduce the main part of your result Theorem 1.1 from Corollary 2.5 of \cite{baldoma2020} except the smoothness at 0 because in higher dimension, in general, the manifold is not smooth at the origin. However it does provide smoothness in $(0,t_0)$.

\textbf{Final remark 2} We are aware of the applications of our results in Celestial Mechanics and Chemistry. We are really happy to hear that there are also applications in Economics. 

\bibliography{references}

\end{document}